\documentclass[12pt]{article}
\usepackage{tikz}
\usepackage{scalefnt}
\usepackage{multicol}
\usepackage{enumerate}

\usepackage[latin1]{inputenc}
\usepackage[T1]{fontenc}
\usepackage{amsmath,amssymb,amsthm}
\usepackage{graphicx}
\usepackage{mathptmx}
\usepackage{enumitem}
\usepackage{array}
\usepackage{dsfont}
\usepackage{txfonts}
\usepackage{color}
\usepackage{tikz}
\usepackage{tikz-3dplot}
\usepackage{caption}
\usepackage{subcaption}
\usepackage{float}
\usetikzlibrary{calc,fadings,decorations.pathreplacing}
\usepackage{verbatim}
\usepackage{cite}
\usepackage{txfonts}
\usepackage{color}

\newtheorem{thm}{Theorem}[section]
\newtheorem{cor}[thm]{Corollary}
\newtheorem{notn}[thm]{Notation}
\newtheorem{lem}[thm]{Lemma}

\newtheorem{defn}[thm]{Definition}
\newtheorem{rem}[thm]{Remark}

\def\R{\mathbb{R}}

\def\S{\mathbb{S}^2_r}

\def\N{\mathbb{N}}
\def\R{\mathbb{R}}

\pgfmathsetmacro{\gr}{(1+sqrt(5))/2}
\pgfmathsetmacro{\igr}{2/(1+sqrt(5))}
\pgfmathsetmacro{\qtwo}{(sqrt(2)}
\pgfmathsetmacro{\qthree}{(sqrt(3)}
\pgfmathsetmacro{\qd}{(1.5*sqrt(2)}

\newcommand\pgfmathsinandcos[3]{%
  \pgfmathsetmacro#1{sin(#3)}%
  \pgfmathsetmacro#2{cos(#3)}%
}
\newcommand\LongitudePlane[3][current plane]{%
  \pgfmathsinandcos\sinEl\cosEl{#2} 
  \pgfmathsinandcos\sint\cost{#3} 
  \tikzset{#1/.estyle={cm={\cost,\sint*\sinEl,0,\cosEl,(0,0)}}}
}
\newcommand\LatitudePlane[3][current plane]{%
  \pgfmathsinandcos\sinEl\cosEl{#2} 
  \pgfmathsinandcos\sint\cost{#3} 
  \pgfmathsetmacro\yshift{\cosEl*\sint}
  \tikzset{#1/.estyle={cm={\cost,0,0,\cost*\sinEl,(0,\yshift)}}} %
}
\newcommand\DrawLongitudeCircle[2][1]{
  \LongitudePlane{\angEl}{#2}
  \tikzset{current plane/.prefix style={scale=#1}}
  \pgfmathsetmacro\angVis{atan(sin(#2)*cos(\angEl)/sin(\angEl))} %
  \draw[current plane] (\angVis:1) arc (\angVis:\angVis+180:1);
  \draw[current plane,dashed] (\angVis-180:1) arc (\angVis-180:\angVis:1);
}
\newcommand\DrawLatitudeCircle[2][1]{
  \LatitudePlane{\angEl}{#2}
  \tikzset{current plane/.prefix style={scale=#1}}
  \pgfmathsetmacro\sinVis{sin(#2)/cos(#2)*sin(\angEl)/cos(\angEl)}
  \pgfmathsetmacro\angVis{asin(min(1,max(\sinVis,-1)))}
  \draw[current plane] (\angVis:1) arc (\angVis:-\angVis-180:1);
  \draw[current plane,dashed] (180-\angVis:1) arc (180-\angVis:\angVis:1);
}

\tikzset{%
  >=latex, 
  inner sep=0pt,%
  outer sep=2pt,%
  mark coordinate/.style={inner sep=0pt,outer sep=0pt,minimum size=3pt,
    fill=black,circle}%
    }

\begin{document}

{\bf \Large \centering  MEASURABLE COLORINGS OF $\S$\\}

\bigskip

{ \centering \Large \bf Greg Malen\\}
{ \centering \large Ohio State University\\}
{ \centering \large malen@math.ohio-state.edu\\}

\bigskip


\section{Introduction}
For a metric space ($M,d$) we consider the infinite graph whose vertices are all points in $M$, and where two points $v,u$ are adjacent if $d(v,u)=1$. Recall that the chromatic number $\chiup(G)$ of a graph $G$ is the minimum $k$ such that there exists a function $\varphi:G\rightarrow\lbrace1,2,\ldots,k\rbrace$ with $\varphi(v)\neq\varphi(u)$ whenever $v$ and $u$ are adjacent. $\chiup(M)$ was classically studied in $\R^2$, where the best current bounds, $4\leq\chiup\left(\R^2\right)\leq 7$, have remained unimproved upon for nearly the problem's entire history. In 1981 Falconer introduced the measurable chromatic number $\chiup_m$, which additionally requires that $\varphi^{-1}(i)$ is a measurable set in $M$ for all $i$. While lower bounds for $\chiup$ usually come from finite graphs, Falconer showed that in this more restrictive setting, $\chiup_m(\R^2)\geq5$ ~\cite{Falconer}. Here we adapt Falconer's techniques to study $\chiup_m$ on $\S$, the 2-sphere of radius $r$.\\
\\
Erd\H{o}s first looked at $\chiup(\S)$, using chordal distances in the ambient space $\R^3$, conjecturing that $\displaystyle \lim_{r\rightarrow\infty}\chiup\left(\S\right)=\chiup\left(\R^2\right)$, and that $\chiup\left(\mathbb{S}^2_{1/\sqrt{3}}\right) = 3$. The latter conjecture was disproved by Simmons, however, who showed that in fact $\chiup(\S) \geq 4$ for all $r\geq\frac{1}{\sqrt{3}}$~\cite{Simm76}. Simmons went on to conjecture that $\chiup\left(\mathbb{S}^2_{1/\sqrt{3}}\right) = 4$, and $\chiup(\S) \geq 4$ for $r > \frac{1}{2}$. In this paper we show that arguments used by Falconer and Croft in the plane can be adapted to the spherical setting to prove that $\chiup_m\left(\S\right)$ is at least 5 for all but a countable set of $r > \frac{1}{\sqrt{3}}$. Furthermore, we examine the monotonicity of $\chiup_m$ as a function of simply connected 2-dimensional surfaces of constant curvature.\\
\\
\noindent Current results for $\chiup\left(\S\right)$ were mostly chronicled by Simmons in ~\cite{Simm74}. The lower bound of 4 for all $r\geq\frac{1}{\sqrt{3}}$ is still the best known. Upper bounds of 4 for $r\leq\frac{\sqrt{3-\sqrt{3}}}{2}$, 5 for $r\leq\frac{1}{\sqrt{3}}$, and 6 for $r\leq\frac{\sqrt{3}}{2}$ may be obtained by, respectively: a tetrahedron, a cap and four triangular stripes, and a dodecahedron. The only value of $r > \frac{1}{2}$ for which an exact answer can be given is $r=\frac{1}{\sqrt{2}}$. At this radius Simmons notes that a curved octahedron yields a 4-coloring that assigns antipodal regions the same color ~\cite{Simm76}. The regions in this coloring are all measurable and there is a lower bound of 4, hence $\chiup\left(\mathbb{S}^2_{1/\sqrt{2}}\right)=\chiup_m\left(\mathbb{S}^2_{1/\sqrt{2}}\right)=4$.\\

\begin{figure}[H]
   \begin{subfigure}[t]{.3\textwidth}
        \begin{tikzpicture} 

\def\R{1.41} 
\def\angEl{15} 
\def\angPhi{-30}
\filldraw[ball color=white] (0,0) circle (\R);
\DrawLatitudeCircle[\R]{\angPhi}

\pgfmathsetmacro\H{\R*cos(\angEl)};
\LatitudePlane[qzplane]{\angEl}{\angPhi};
\LongitudePlane[azplane]{\angEl}{-150};
\LongitudePlane[bzplane]{\angEl}{120};
\tdplotdrawarc[azplane]{(0,0,0)}{\R}{118}{-30}{}{};
\tdplotdrawarc[azplane, dashed]{(0,0,0)}{\R}{118}{210}{}{};
\tdplotdrawarc[bzplane]{(0,0,0)}{\R}{210}{72}{}{};
\tdplotdrawarc[bzplane, dashed]{(0,0,0)}{\R}{72}{-30}{}{};
\end{tikzpicture}
\caption{$r=\frac{1}{\sqrt{3}}; \ \chiup\leq5$}
    \end{subfigure}
\hspace{-.4cm}
    \begin{subfigure}[t]{.3\textwidth}
    	\begin{tikzpicture}
	\def\R{1.73} 
\def\angEl{15} 
\def\angPhi{-30}
\filldraw[ball color=white] (0,0) circle (\R);
\DrawLatitudeCircle[\R]{0}

\pgfmathsetmacro\H{\R*cos(\angEl)}
\DrawLongitudeCircle[\R]{-150}
\DrawLongitudeCircle[\R]{120}
\end{tikzpicture}
\caption{$r=\frac{1}{\sqrt{2}}; \ \chiup=4$}
\end{subfigure}
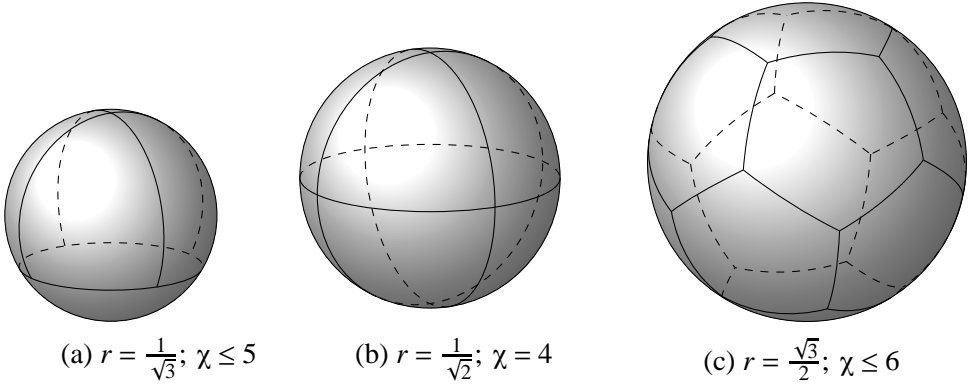
\hspace{.3cm}
    \begin{subfigure}[t]{.3\textwidth}
    \begin{tikzpicture} 

\def\R{2.12} 
\def\angEl{15} 
\def\angAlpha{46.78}
\def\angZeta{-86}
\def\angGamma{58}
\def\angPhi{-52.62}
\def\angBeta{52.62}
\def\angLatPts{10.81}
\filldraw[ball color=white] (0,0) circle (\R);

\pgfmathsetmacro\H{\R*cos(\angEl)}

\LatitudePlane[qzplane]{\angEl}{\angPhi};
\LongitudePlane[azplane]{\angEl}{-42};
\LongitudePlane[bzplane]{\angEl}{30};
\LongitudePlane[czplane]{\angEl}{102};
\LongitudePlane[dzplane]{\angEl}{174};
\LongitudePlane[ezplane]{\angEl}{-114};
\LongitudePlane[midAplane]{\angEl}{-132}
\LongitudePlane[midBplane]{\angEl}{120}
\LongitudePlane[midCplane]{\angEl}{12}
\LongitudePlane[midDplane]{\angEl}{-96}
\LongitudePlane[midEplane]{\angEl}{-24}

\tdplotdrawarc[azplane]{(0,0,0)}{\R}{10.81}{52.62}{}{};
\tdplotdrawarc[azplane, dashed]{(0,0,0)}{\R}{190.81}{232.62}{}{};
\tdplotdrawarc[bzplane, dashed]{(0,0,0)}{\R}{10.81}{52.62}{}{};
\tdplotdrawarc[bzplane]{(0,0,0)}{\R}{190.81}{232.62}{}{};
\tdplotdrawarc[czplane, dashed]{(0,0,0)}{\R}{10.81}{52.62}{}{};
\tdplotdrawarc[czplane]{(0,0,0)}{\R}{190.81}{232.62}{}{};
\tdplotdrawarc[dzplane, dashed]{(0,0,0)}{\R}{10.81}{25}{}{};
\tdplotdrawarc[dzplane]{(0,0,0)}{\R}{25}{52.62}{}{};
\tdplotdrawarc[dzplane]{(0,0,0)}{\R}{190.81}{202}{}{};
\tdplotdrawarc[dzplane, dashed]{(0,0,0)}{\R}{202}{232.62}{}{};
\tdplotdrawarc[ezplane]{(0,0,0)}{\R}{10.81}{52.62}{}{};
\tdplotdrawarc[ezplane, dashed]{(0,0,0)}{\R}{190.81}{232.62}{}{};

\path[midAplane] (0:\R) coordinate (A1);
\path[azplane] (212:\R) coordinate (A2);
\path[azplane] (122:\R) coordinate (A3);

 \begin{scope}[ x={(A1)}, y={(A2)}, z={(A3)}]       
          \draw[black] ( -20.905:1) arc (-20.905:20.905:1) ;
          \draw[black, dashed] ( 159.905:1) arc (159.905:200.905:1) ;  
 \end{scope} 

\path[midBplane] (0:\R) coordinate (B1);
\path[bzplane] (212:\R) coordinate (B2);
\path[bzplane] (122:\R) coordinate (B3);

 \begin{scope}[ x={(B1)}, y={(B2)}, z={(B3)}]       
          \draw[black, dashed] ( -20.905:1) arc (-20.905:20.905:1) ;
          \draw[black] ( 159:1) arc (159:200.905:1) ;  
 \end{scope}
 
\path[midCplane] (0:\R) coordinate (C1);
\path[czplane] (212:\R) coordinate (C2);
\path[czplane] (122:\R) coordinate (C3);

 \begin{scope}[ x={(C1)}, y={(C2)}, z={(C3)}]       
          \draw[black, dashed] ( -20.905:1) arc (-20.905:12:1) ;
          \draw[black] ( 12:1) arc (12:20.905:1) ;
          \draw[black] ( 159:1) arc (159:192:1) ;  
          \draw[black, dashed] ( 192:1) arc (192:200.905:1) ; 
 \end{scope}

\path[midDplane] (0:\R) coordinate (D1);
\path[dzplane] (212:\R) coordinate (D2);
\path[dzplane] (122:\R) coordinate (D3);

 \begin{scope}[ x={(D1)}, y={(D2)}, z={(D3)}]       
          \draw[black] ( -20.905:1) arc (-20.905:20.905:1) ;
          \draw[black, dashed] ( 159.905:1) arc (159.905:200.905:1) ;  
 \end{scope} 
 
\path[midEplane] (0:\R) coordinate (E1);
\path[ezplane] (212:\R) coordinate (E2);
\path[ezplane] (122:\R) coordinate (E3);

 \begin{scope}[ x={(E1)}, y={(E2)}, z={(E3)}]       
          \draw[black] ( -20.905:1) arc (-20.905:20.905:1) ;
          \draw[black, dashed] ( 159.905:1) arc (159.905:200.905:1) ;  
 \end{scope}
 
\path[azplane] (211.75:\R) coordinate (a2);
\path[azplane] (121.75:\R) coordinate (a3);
\path[bzplane] (211.75:\R) coordinate (b2);
\path[bzplane] (121.75:\R) coordinate (b3);
\path[czplane] (211.75:\R) coordinate (c2);
\path[czplane] (121.75:\R) coordinate (c3);
\path[dzplane] (211.75:\R) coordinate (d2);
\path[dzplane] (121.75:\R) coordinate (d3);
\path[ezplane] (211.75:\R) coordinate (e2);
\path[ezplane] (121.75:\R) coordinate (e3);

  \begin{scope}[ x={(A1)}, y={(a3)}, z={(a2)}]       
          \draw[black] ( 69:1) arc (69:80:1) ;
          \draw[black, dashed] ( 80:1) arc (80:110.905:1) ;
          \draw[black] ( -69:1) arc (-69:-98:1) ;
          \draw[black, dashed] ( -98:1) arc (-98:-110.905:1) ;  
 \end{scope}
 
   \begin{scope}[ x={(B1)}, y={(b3)}, z={(b2)}]       
          \draw[black] ( 69:1) arc (69:110.905:1) ;
          \draw[black, dashed] ( -69:1) arc (-69:-110.905:1) ;  
 \end{scope}

  \begin{scope}[ x={(C1)}, y={(c3)}, z={(c2)}]       
          \draw[black] ( 69:1) arc (69:110.905:1) ;
          \draw[black, dashed] ( -69:1) arc (-69:-110.905:1) ;  
 \end{scope}
 
   \begin{scope}[ x={(D1)}, y={(d3)}, z={(d2)}]       
          \draw[black] ( 69:1) arc (69:110.905:1) ;
          \draw[black, dashed] ( -69:1) arc (-69:-110.905:1) ;  
 \end{scope}
 
   \begin{scope}[ x={(E1)}, y={(e3)}, z={(e2)}]       
          \draw[black, dashed] ( 69:1) arc (69:110.905:1) ;
          \draw[black] ( -69:1) arc (-69:-110.905:1) ;  
 \end{scope}

\end{tikzpicture}
\caption{$r=\frac{\sqrt{3}}{2}; \ \chiup\leq6$}
    \end{subfigure}
\end{figure}

\noindent For any point $x\in\S$, the set of points unit distance from $x$ forms a circular cross-section with radius $R=\frac{\sqrt{4r^2-1}}{2r}$. On this circle, the central angle subtended by two points which are themselves unit distance apart is $\theta=2\sin^{-1}\left(\frac{1}{2R}\right) = 2\sin^{-1}\left(\frac{r}{\sqrt{4r^2-1}}\right)$.
\begin{thm}
For $r > \frac{1}{\sqrt{3}}$, let $\theta = 2\sin^{-1}\left(\frac{r}{\sqrt{4r^2-1}}\right) = c\cdot\pi$. If $c$ is irrational, or if $c=\frac{p}{q}\cdot\pi$ with $p,q\in\N, \ (p,q)=1$ and $p$ is even, then $\chiup_m(\S)\geq5$.
\end{thm}
\vskip .5cm
\noindent The set of $r$ values not covered by this theorem, namely when $p$ is odd, is only a countable set. These values, being outside the purview of the methods used in this paper, may or may not provide exceptions to the theorem. The octahedral coloring for $r=\frac{1}{\sqrt{2}}$ yields the only known example of such an exception, with $\theta=\frac{\pi}{2}$.

\vskip 1cm
\section{Measure Theoretic Lemmas}
\noindent First we introduce some measure theoretic notations and definitions.

\begin{notn}
Let $\mu$ and $l$ be the 2-dimensional and 1-dimensional spherical measures on $\mathbb{S}_r^2$ respectively, and let $\alpha$ be 1-dimensional angular measure.
\end{notn}

\begin{defn}
For $E$ a measurable set, the density at a point $x\in E$ is $\delta_E(x)\coloneqq\displaystyle\lim_{R\rightarrow0}\frac{\mu(B_{R}(x)\cap E)}{\mu(B_{R}(x))}$.
\end{defn}

\begin{defn}
For $E$ a measurable set, define its metrical boundary $\partial E\coloneqq\lbrace x\in\S : \delta_E(x)$ does not exist, or \ $0 < \delta_E(x) < 1\rbrace$, and the essential part $\tilde{E} \coloneqq\lbrace x\in\S : \delta_E(x)=1\rbrace$.
\end{defn}

\begin{rem}
Let $E$ be a measurable subset of $\S$ with $\mu(E) > 0$ and $\mu(\S\setminus E) > 0$.  Then $\partial E\neq\varnothing, \ \mu(\partial E) = 0,$ and $\tilde{E}$ is a Borel set.~\cite{Falconer}
\end{rem}

\vskip.3cm
\noindent The following lemma is modeled off of Croft's argument for lattices in the plane~\cite{Croft}. It ensures for a particular arrangement of points $L$, and any set $E$ with zero measure, that $|E\cap L| = 1$ for a natural set of rigid motions of $\S$. Here we let $L$ be 2 points at unit chordal distance in $\S$, and the rigid motions are rotations of $\S$ about points, where rotating about a point $x$ means precisely rotating the sphere about the diameter that goes through $x$ and its antipodal point.\\

\begin{lem}
Fix $r > \frac{1}{2}$ and let $E\subset\S$ with $\mu(E) = 0$. Let $L = \lbrace x_1,x_2\rbrace\subset\S$ be two points at unit chordal distance.  Then for a given $x_i$ there is an orientation of $L$ in $\S$ with $x_i\in E\cap L$, such that for almost all rotations $\rho$ about $x_i, \ \rho(E)\cap L=\lbrace x_i\rbrace$.
\end{lem}

\noindent{\bf Proof of Lemma. } Assume the statement is false. So for every orientation of $L$ with a given point $x_i\in E\cap L$ there is a set $C$ of rotations about $x_i$ such that $\alpha(C) > 0$ with $\rho(E)\cap L=L$ for all $\rho\in C$.\\
\\
Without loss of generality, place $L$ on $\S$ such that $x_1 \in E\cap L$.  Note that $r > \frac{1}{2}$ ensures that points at unit chordal distance exist on $\S$. By assumption, there is a set $C$ of rotations about $x_1$ with $\alpha(C) > 0$, for which $x_2\in\rho(E)\cap L$ \ for all $\rho\in C$.  And since $r > \frac{1}{2}, \ x_1 \text{ and } x_2$ are not antipodal, hence $x_2$ traces out a circular cross section of $\S$ via rotations about $x_1$. Then $\alpha(C) > 0$ if and only if there is a set of points $D\subset E$ on this cross section such that $l(D) > 0$.\\
\\
Now reorient $L$ so that $x_1\in D$. Again, for $x_1$ starting at any point $y\in D$ there is a set $C_y$ of rotations about $y$, with $\alpha(C_y) > 0$, for which $x_2\in\rho(E)\cap L$ \ for all $\rho\in C_y$. As before, this produces a set of points $D_y\subset E$ on the cross section of points chordal distance 1 from $y$, such that $l(D_y) > 0$.\\
\\
Now consider the integral $$\int_{D}\left(\int_{D_y}\mathds{1}_E \ \text{d}(l)\right)\text{d}(l)=\mu\left(\bigcup_{y\in D} D_y\right)$$ Recall that the integral of a measurable, non-negative function is zero only if the function is zero on all but a set of measure zero. Here for every $y\in D$ $$\int_{D_y}\mathds{1}_E \ \text{d}(l)=l(D_y) > 0$$ Therefore \ $\displaystyle\mu(E)\geq\mu\left(\bigcup_{y\in D} D_y\right) > 0$, which is a contradiction. $_{\square}$\\
\\
\\
It should be noted that this lemma does not hold in general for any discrete set $L\subset\S$. In particular, the lemma is false for arrangements where $E$ and $L$ are both such that $x\in E, L$ respectively if and only if its antipodal point is also in $E, L$. This contrasts the planar case, where Croft showed the lemma holds for any lattice $L$ and any set $E$ with zero measure~\cite{Croft}. The following corollary to Lemma 2.5 and pair of lemmas are due to Falconer~\cite{Falconer},~\cite{Soifer}, however the proof of the first lemma must be adapted to fit the geometry of the sphere.\\
\begin{cor}
Let $\mathbb{S}_r^2$ be covered by disjoint measurable sets $S_1,\ldots,S_m$, and let $E=\S\setminus\displaystyle\bigcup_{i=1}^m\tilde{S}_i$.  Then $\mu(E)=0$, and there is an orientation of the sphere with $x_1\in E, \ x_2\in\displaystyle\bigcup_{i=1}^m\tilde{S}_i$, \ such that for almost all rotations $\rho$ about $x_1, \ x_2$ remains in \ $\displaystyle\bigcup_{i=1}^m\tilde{S}_i$.
\end{cor}

\vskip.4cm
\noindent Observe that for $x\in E, \ \delta_i(x)$ may be 0 for at most two of the $S_i$. If it were 0 for three of the four sets, then the density with respect to the last set would be 1, and $x$ would not be in $E$. Hence if $x\in E$ then $x\in\displaystyle\bigcap_{i=1}^{k}\partial S_{n_i}$, for some $\lbrace n_1,\ldots,n_k\rbrace\subset[m]\coloneqq\lbrace1,2,\ldots,m\rbrace, \ 2\leq k\leq m$.

\vskip1cm
\begin{lem}
Let $\mathbb{S}_r^2$ be covered by disjoint measurable sets $S_1,\ldots,S_m$, and let $L$ be in the orientation given by the corollary. Let $x_1\in\displaystyle\bigcap_{i=1}^{k}\partial S_{n_i}$, \ with $\lbrace n_1,\ldots,n_k\rbrace\subset[m], \ 2\leq k\leq m$, and let $x_2\in\tilde{S}_{\sigma(2)}$ \ for some $\sigma(2)\in[m]$. Then for every $i\in[k]$ there exists a rigid motion $\rho_i$ of $\S$ which moves $x_1$ into $S_{n_i}$ and $x_2$ into $S_{\sigma(2)}$. 
\end{lem}
\vskip .5cm
\noindent{\bf Proof.} Without loss of generality, let $x_1$ start at the north pole. For a given $R > 0$ we can associate each point in the ball of radius $R$ about $x_1, \ B_R(x_1)$, to a unique rigid motion which rotates $x_1$ directly along the spherical geodesic to this point. All of these rigid motions will be rotations about points on the equator, as they rotate $x_1$ along meridians of $\S$.  Call this set of rigid motions $P_R$. As $\mu(B_R(x))$ depends only on $R$, define $A_R=\mu(B_R(x))$ for any $x\in\S$. Then by definition of density, there exists some $\epsilon > 0$ small enough that for some arbitrarily small $R$, the following two inequalities both hold.\\
\\
\hspace*{3.4cm}$\displaystyle \epsilon \ < \ \ \ \frac{\mu\left(B_R(x_1)\cap S_{n_i}\right)}{A_R} \ \ \ < \ 1-\epsilon$\\
\\
\hspace*{2.8cm}$\displaystyle 1-\frac{\epsilon}{2} \ \leq \ \frac{\mu\left(B_R(x_2)\cap S_{\sigma(2)}\right)}{A_R} \ \leq \ 1$\\
\\
Note that $d(\rho(x_2),x_2)\leq d(\rho(x_1),x_1)$ for all $\rho\in P_R$, with equality only holding when $\rho(x_1)$ and $x_2$ are on the same meridian. Then for $P_2\coloneqq\lbrace\rho(x_2):\rho\in P_R\rbrace\subset B_R(x_2)$, $$\frac{\mu\left(P_2\setminus\left(P_2\cap S_{\sigma(2)}\right)\right)}{A_R}\leq\frac{\mu\left(B_R(x_2)\setminus\left(B_R(x_2)\cap S_{\sigma(2)}\right)\right)}{A_R} < \frac{\epsilon}{2}$$ Define $P_{i,j}\coloneqq\lbrace x\in B_R(x_1)\cap S_{n_i}: \rho(x_2)\in S_j \text{ when } x=\rho(x_1)\rbrace$.\\
\\
$$\displaystyle\frac{\mu\left(P_{i,\sigma(2)}\right)}{A_R}\geq \ \frac{\mu\left(B_R(x_1)\cap S_{n_i}\right)}{A_R} - \frac{\mu\left(P_2\setminus(P_2\cap S_{\sigma(2)})\right)}{A_R}$$
\hspace*{3.3cm}$\displaystyle \ > \ \epsilon-\frac{\epsilon}{2}$\\
\\
\hspace*{3.3cm} $> \ 0$
\\
\\ 
In particular for every $i, \ P_{i,\sigma(2)}\neq\emptyset. \ _{\square}$\\

\begin{lem}
Let $R$ be such that $R > \frac{1}{2}$ and $\theta=2\sin^{-1}\left(\frac{1}{2R}\right)$ is an irrational multiple of $\pi$. Let $S_1,S_2$ be measurable subsets of $\mathbb{R}^2$. Suppose almost all points on a circle of radius $R$ are in either $\tilde{S}_1$ or $\tilde{S}_2$.  Then at least one of $S_1,S_2$ contains a pair of points that are unit distance apart on the circle.
\end{lem}

\noindent For a proof see Falconer, Lemma 4~\cite{Falconer}.
\section{Proof of Theorem.}

\noindent Fix $r > \frac{1}{\sqrt{3}}$. Let $\theta=2\sin^{-1}\left(\frac{r}{\sqrt{4r^2-1}}\right) = c\cdot\pi$, and suppose that $c$ is irrational. Suppose there is a permissible measurable 4-coloring of $\S$, with disjoint color classes $S_1,S_2,S_3,S_4$.  Let $x_1,x_2\in\S$ be a pair of points at unit chordal distance from each other.  By Corollary 2.6, there is an orientation of the sphere in which $x_1\in\S\setminus\displaystyle\bigcup_{i=1}^4\tilde{S}_i, \ \text{and} \ x_2\in\bigcup_{i=1}^4\tilde{S}_i$ \ for almost all rotations about $x_1$.  Without loss of generality assume $x_1\in\partial S_1\cap\partial S_2$. So $x_2\notin\tilde{S}_1\cup\tilde{S}_2$, since Lemma 2.7 gives rigid motions which would place it in the same color class as $x_1$.  Hence $x_2\in\tilde{S}_3\cup\tilde{S}_4$ for almost all rotations about $x_1$, through which $x_2$ traces out a cross-section of $\S$ which is a circle in $\mathbb{R}^2$. By Lemma 2.8 and the assumption that $\theta$ is an irrational multiple of $\pi$, either $S_3$ or $S_4$ must contain a pair of points on this circle unit distance apart, which is a contradiction.\\
\\
Now suppose that $c=\frac{p}{q}$ for $p,q\in\N$ with $p$ even and $(p,q)=1$, hence $q$ odd. Let $N$ be the north pole of the sphere and consider the circular cross section of points unit distance from $N$. Let $L=\lbrace N,y_1,\ldots,y_{q}\rbrace$, where the $y_j$ are evenly spaced along the cross section at rotations of $\frac{2\pi}{q}$. Since $(p,q)=1$, the smallest multiple of $\theta=\frac{p\pi}{q}$ which is an integer multiple of $\pi$ is $q\cdot\theta=p\pi$. And $p$ is even so this is a multiple of $2\pi$. Hence rotating by $\theta$ forms a unit-distance odd cycle among the $q$ points evenly spaced at rotations of $\frac{2\pi}{q}$ in some order. Since the $y_j$ trace out the same cross section as $x_2$ did, Lemma 2.5 still holds for $N\in E\cap L$, as the pigeonhole principle guarantees that one of the $y_j$ would give an arc of positive measure at each step. So there is an orientation with $\rho(E)\cap L=\lbrace N\rbrace$ for almost all rotations about $N$. Then for $y_j\in\tilde{S}_{\sigma(j)}$ for each $j$, Lemma 2.7 must be adjusted to find rotations $\rho_1$ and $\rho_2$ which move $x_1$ into $S_1$ and $S_2$ respectively, while simultaneously moving $y_j$ into $S_{\sigma(j)}$ for every $j$. We can define $P_j\coloneqq \lbrace\rho(y_j):\rho\in P_R\rbrace$ for each $j$, and choose $\frac{\epsilon}{q+1}$ instead of $\frac{\epsilon}{2}$. The subtraction in the proof then becomes $$\frac{\mu\left(B_R(x_1)\cap S_{n_i}\right)}{A_R} - \sum_{j=1}^q\frac{\mu\left(P_i\setminus\left(P_j\cap S_{\sigma(j)}\right)\right)}{A_R} > \ \epsilon-\displaystyle\sum_{j=1}^{q}\frac{\epsilon}{q+1}$$
This is greater than 0, so the desired $\rho_1$ and $\rho_2$ exist. As before, this yields that if there is a permissible 4-coloring, then each $y_j$ is in $\tilde{S}_3\cup \tilde{S}_4$, and are moved into $S_3\cup S_4$ by $\rho_1$ and $\rho_2$. But this is a contradiction since the $y_j$ form an odd cycle. Hence there is no permissible 4-coloring. $_{\square}$\\

\section{Covering More Ground}
\noindent Despite the exception when $\theta=\frac{\pi}{2}$, it may be possible to cover a large set of the cases for which $\theta=\frac{p\pi}{q}$ with $p$ odd by looking at other unit distance configurations. The unit distance diamond, $K_4$ minus an edge, is one such example.

\begin{thm}
Let $L=\lbrace x_1,x_2,x_3,x_4\rbrace$ be a unit distance graph with all possible edges except for $x_1x_4$. With $x_1$ at the north pole, let $D$ be the radius of the horizontal cross section traversed by $x_4$. Then $\ \chiup_m(\S)\geq5$ when $D > \frac{1}{2} \text{ and }\beta=2\sin^{-1}\left(\frac{1}{2D}\right)$ is an irrational multiple of $\pi$ or a rational multiple of $\pi$ with an even numerator. 
\end{thm}

\noindent Lemma's 2.5 and 2.7 can be adapted to fit this choice of $L$ as in the proof of the theorem, so long as $x_1$ and $x_4$ are not antipodal. But that is the exceptional case $r=\frac{1}{\sqrt{2}}$. Then without loss of generality $x_1\in\partial S_1\cap\partial S_2$, and so $x_2,x_3\in\tilde{S}_3\cup\tilde{S}_4$. Thus $x_4$ must be in $\tilde{S}_1\cup\tilde{S}_2$ for almost all rotations about $x_1$, and the proof follows as before.  $D$ was initially computed by Simmons ~\cite{Simm76} to be $$D=\left|\frac{2r(2r^2-1)\sqrt{3r^2-1}}{4r^2-1}\right|$$ So $D > \frac{1}{2}$ for $r\in(0.586158,0.627745)\cup(0.819417,\infty)$, where these are decimal approximations of the solutions to $D=\frac{1}{2}$. $\beta$ can be written in terms of $\theta$ as $$\beta=2\sin^{-1}\left(\frac{\tan\left(\frac{\theta}{2}\right)}{4\cos(\theta)}\right)$$ However there is not a simple way to check when $\sin^{-1}(x)$ is a rational multiple of $\pi$ in general, making it difficult to decide when the $r$ values covered by these two theorems intersect.
\\
\section{Further Questions and Implications}

An intriguing corollary of Theorem 1.1 is that $\chiup_m(\S)$ is not monotonic as a function of $r$, which gives rise to further questions as to the nature of $\chiup$ and $\chiup_m$ as functions of curvature. For negatively curved spaces, using a varying distance $d$ to define the edge set instead of varying the curvature, Kloeckner showed the following bounds for $\chiup(\mathbb{H}_d^2)$ ~\cite{Kloeckner}. 

\begin{thm}
For all $d\geq3\ln(3), \ \chiup(\mathbb{H}_d^2)\leq4\left\lceil\frac{d}{\ln(3)}\right\rceil+4$. And for $d\leq3\ln(\frac{3}{2}), \ \chiup(\mathbb{H}_d^2)\leq12$.
\end{thm}

\noindent The proof uses a checkerboard lattice to cover $\mathbb{H}_d^2$. The more natural regular n-gon tilings of $\mathbb{H}^2$ with three n-gons meeting at a vertex, as with hexagons in the plane, surprisingly do much worse. In this setting we ask whether there is a lower bound which also grows with $d$, and whether $\displaystyle \lim_{c\rightarrow0}\chiup\left(\mathbb{H}_c^2\right)=\chiup\left(\R^2\right)$, where $c$ is the fixed constant curvature. Furthermore, in the spherical setting we ask if there is any $r$ such that $\chiup_m(\S)=5$, and if there is any $r$ besides $r=\frac{1}{\sqrt{2}}$ where $\chiup_m(\S)=4$.
\bibliographystyle{plain}
\bibliography{Master}
\end{document}